\documentclass[11pt,leqno]{article}
\usepackage{amsmath,amstext}

\date{11 June 1997}

\parfillskip=\parindent plus1fil %prevents last line of paragraph being full

\renewcommand{\phi}{\varphi}

\newcommand{\ip}[1]{\langle #1 \rangle}

\newcommand{\T}{\mathcal{T}}
\newcommand{\B}{\mathcal{B}}
\newcommand{\M}{\mathcal{M}}
\newcommand{\U}{\mathcal{U}}
\newcommand{\CT}{\mathcal{C_T}}
\newcommand{\CM}{C(\M)}
\newcommand{\wt}{\widetilde}    
\newcommand{\tl}{\tilde}
\newcommand{\til}{\tl{\ }}
\newcommand{\hop}{\text{\textsl{HOP}}}

\DeclareMathOperator{\re}{Re}

\newcommand{\bibref}[1]{[\ref{#1}]}

\renewcommand{\emph}{\textsl}

\newcommand{\proof}{\textsc{Proof:} \hspace{.4em minus .2em}}

\newcommand{\qed}{\parfillskip=0pt plus1fil\nolinebreak\hfill$\rule{1ex}{1ex}$%
\par\addvspace{12pt plus 3pt minus 3pt}\parfillskip=\parindent plus1fil}

\title{\sffamily\bfseries\Large The Berezin Transform on the Toeplitz Algebra}
\author{\sffamily\bfseries Sheldon Axler \and \sffamily\bfseries Dechao Zheng}

\newenvironment{keeptogether}{\pagebreak[0]\samepage}{}

\numberwithin{equation}{section}

\newtheorem{Theorem}[equation]{\sffamily\bfseries Theorem}
\newtheorem{Proposition}[equation]{\sffamily\bfseries Proposition}
\newtheorem{Lemma}[equation]{\sffamily\bfseries Lemma}
\newtheorem{Corollary}[equation]{\sffamily\bfseries Corollary}
\newtheorem{Conjecture}[equation]{\sffamily\bfseries Conjecture}

\newenvironment{theorem}{\begin{keeptogether}\begin{Theorem}\slshape}%
  {\end{Theorem}\end{keeptogether}}
\newenvironment{proposition}{\begin{keeptogether}\begin{Proposition}\slshape}%
  {\end{Proposition}\end{keeptogether}}
\newenvironment{lemma}{\begin{keeptogether}\begin{Lemma}\slshape}%
  {\end{Lemma}\end{keeptogether}}
\newenvironment{corollary}{\begin{keeptogether}\begin{Corollary}\slshape}%
  {\end{Corollary}\end{keeptogether}}
  {\end{Conjecture}\end{keeptogether}}

\newcounter{referencec}

\newcounter{subtheoremc}
\newcounter{subsubtheoremc}

\newenvironment{subtheorem}{\begin{list}{{\normalfont
(\alph{subtheoremc})}\hfill}{\usecounter{subtheoremc}
\setlength{\labelwidth}{.4in}
\setlength{\labelsep}{0pt}
\setlength{\leftmargin}{.4in}
}}{\end{list}}

\newenvironment{subsubtheorem}{\begin{list}{{\normalfont
(\roman{subsubtheoremc})}\hfill}{\usecounter{subsubtheoremc}
\setlength{\labelwidth}{.3in}
\setlength{\labelsep}{0pt}
\setlength{\leftmargin}{.3in}
}}{\end{list}}

\makeatletter

\renewcommand{\@startsection}[6]{\if@noskipsec \leavevmode
\fi
   \par \@tempskipa #4\relax
   \@afterindenttrue
   \ifdim \@tempskipa <\z@ \@tempskipa -\@tempskipa \@afterindenttrue\fi
   \if@nobreak \everypar{}\else
     \addpenalty{\@secpenalty}\addvspace{\@tempskipa}\fi \@ifstar
     {\@ssect{#3}{#4}{#5}{#6}}{\@dblarg{\@sect{#1}{#2}{#3}{#4}{#5}{#6}}}}

\renewcommand{\section}{\@startsection {section}{1}{\z@}%
                                   {-3.5ex \@plus -1ex \@minus -.2ex}%
                                   {2.3ex \@plus.2ex}%
                                   {\reset@font\large\sffamily\bfseries}}

\makeatother

\renewcommand{\[}{\parfillskip=0pt plus1fil$$}
\renewcommand{\]}{$$\parfillskip=\parindent plus1fil}

\makeatletter
\renewenvironment{equation}{\parfillskip=0pt plus1fil%
\refstepcounter{equation}%
$$%
}{%
\leqno\tagform@{\theequation}%
$$%
\parfillskip=\parindent plus1fil}%

\newcommand{\simpletag}[1]{\def\@currentlabel{#1}\def\theequation{#1}%
  \addtocounter{equation}{-1}}
\makeatother

\begin{document}
\maketitle
\thispagestyle{empty}

\begin{list}{}{\setlength{\leftmargin}{\parindent}
\setlength{\rightmargin}{\parindent}}

\item \small
\textit{Abstract}.  This paper studies the boundary behavior of the Berezin
transform on the $C^*$-algebra generated by the analytic Toeplitz operators on the
Bergman space.
\end{list}

\section{Introduction}

{}Let $dA$ denote Lebesgue area measure on the unit disk $D$, normalized so that the
measure of~$D$ equals $1$. The \emph{Bergman space} $L^{2}_{a}$ is the Hilbert
space consisting of the analytic functions on~$D$ that are also in~$L^{2}(D, dA)$.
For $z \in D$, the \emph{Bergman reproducing kernel} is the function $K_z \in L^2_a$
such that\footnote{Both authors
were partially supported by the National Science Foundation.  The first author also
thanks the Mathematical Sciences Research Institute (Berkeley), for its hospitality
while part of this work was in progress.}
\[
f(z) = \ip{f, K_z}
\]
for every $f \in L^2_a$.  The \emph{normalized Bergman reproducing kernel} $k_z$ is the
function $K_z/\|K_z\|_2$.  Here, as elsewhere in this paper, the norm $\|\ \|_2$ and
the inner product
$\ip{\ \,,\ }$ are taken in the space $L^2(D, dA)$.  The set of bounded operators
on $L^2_a$ is denoted by~$\B(L^2_a)$.

For $S \in \B(L^2_a)$, the \emph{Berezin
transform} of $S$ is the function $\tl{S}$ on $D$ defined by 
\[
\tl{S}(z)=\ip{Sk_{z}, k_{z}}.
\]
Often the behavior of the Berezin transform of an operator provides
important information about the operator.

For $u\in L^\infty(D, dA)$, the \emph{Toeplitz operator} $T_u$ with symbol $u$ is the
operator on $L^2_a$ defined by $T_{u}f = P(uf)$; here $P$ is the orthogonal projection
from $L^{2}(D, dA)$ onto $L_{a}^{2}$.  Note that if $g \in H^\infty$ (the set of
bounded analytic functions on~$D$), then $T_g$ is just the operator of multiplication
by $g$ on~$L^2_a$.

The Berezin transform  $\tl{u}$ of a function $u \in L^\infty(D, dA)$ is defined to
be the Berezin transform of the Toeplitz operator $T_u$.  In other words,
$\tl{u} = \wt{T_u}$.  Because $\ip{T_u k_z, k_z} = \ip{P(uk_z), k_z} =
\ip{uk_z, k_z}$, we
obtain the formula
\begin{equation} \label{49}
\tl{u}(z) = \int_D u(w) |k_z(w)|^2\,dA(w).
\end{equation}%
The Berezin transform of a function in $L^\infty(D, dA)$ often
plays the same important role in the theory of Bergman spaces as the harmonic
extension of a function in $L^\infty(\partial D, d\theta)$ plays in the theory of Hardy
spaces.

The \emph{Toeplitz algebra} $\T$ is the $C^*$-subalgebra of $\B(L^2_a)$ generated by
\mbox{$\{T_g : g \in H^\infty\}$}.  We let $\U$ denote
the \mbox{$C^*$-subalgebra} of $L^\infty(D, dA)$ generated by~$H^\infty$.  As is
well known (see \bibref{Ax2}, Proposition~4.5), $\U$ equals the closed subalgebra
of $L^\infty(D, dA)$ generated by the set of bounded harmonic functions on~$D$. 
Although the map $u \mapsto T_u$ is not multiplicative on $L^\infty(D, dA)$, the
identities ${T_u}^* = T_{\bar{u}}$, $T_u T_g = T_{ug}$, and
$T_{\bar g} T_u = T_{\bar{g}u}$ hold for all $u \in L^\infty$ and all
$g \in H^\infty$.  This implies that $\T$ equals the closed subalgebra of
$\B(L^2_a)$ generated by the Toeplitz operators with bounded harmonic symbol, and
that $\T$ also equals the closed subalgebra of $\B(L^2_a)$ generated by
\mbox{$\{T_u : u \in \U\}$}.  Our goal in this paper is to study the boundary behavior
of the Berezin transforms of the operators in
$\T$ and of the functions in~$\U$.

In Section 2 we study the boundary behavior of Berezin transforms of operators
in~$\T$.  We show (Theorem~\ref{23}) that if $S \in \T$, then
$\tl{S} \in \U$.  Perhaps the main result in this section is Theorem~\ref{22}, which
describes the commutator ideal $\CT$ (the smallest closed, two-sided ideal of
$\T$ containing all operators of the form $RS - SR$, where $R, S \in \T$).  As a
consequence of this result, we show (Corollary~\ref{66}) that $S - T_{\tl{S}}$ is in
the commutator ideal $\CT$ for every $S \in \T$.  Writing $S = T_{\tl{S}}
+ (S - T_{\tl{S}})$, this gives us a canonical way to express the (nondirect) sum $\T =
\{T_u : u \in \U\} + \CT$.  We also prove (Corollary~\ref{77}) that if $S \in \CT$,
then $\tl{S}$ has nontangential limit $0$ at almost every point of~$\partial D$.

In Section~3 we study the boundary behavior of Berezin transforms of functions
in~$\U$.  We prove (Corollary~\ref{42}) that if $u \in \U$, then $\tl{u} - u$ has
nontangential limit $0$ at almost every point of~$\partial D$.  Using similar
techniques, we prove (Corollary~\ref{43}) that if $u \in \U$, then the function
$z \mapsto \|T_{u - u(z)} k_z\|_2$ has nontangential limit $0$ at almost every point
of~$\partial D$.  The main result of this section is Theorem~\ref{36}, which
describes the functions $u \in \U$ such that $\tl{u}(z) - u(z) \to 0$ as $z \to
\partial D$.  As a consequence, we describe (Corollary~\ref{44}) the
operators $S$ that differ from the Toeplitz operator $T_{\tl{S}}$ by a compact
operator, where $S$ is a finite sum of finite products of Toeplitz operators with
symbols in~$\U$.

In Section~4 we use results from the two previous sections to describe
(Theorem~\ref{45}) when the Berezin transform is asymptotically multiplicative on
harmonic functions.  This theorem is then used to characterize the functions $f, g \in
H^\infty$ such that $T_{\bar{f}} T_g - T_g T_{\bar{f}}$ is compact.

We thank Jaros{\l}aw Lech for useful conversations about the Berezin transform.

\section{Boundary Behavior of the Berezin Transform on $\T$}

In this section we will study the boundary behavior of the Berezin transform on
elements of~$\T$.  We will need explicit formulas for the reproducing kernel and the
normalized reproducing kernel.  As is well known,
\[
K_z(w) = \frac{1}{(1 - \bar{z}w)^2}
\]
for $z, w \in D$.  Note that
\[
{\|K_z\|_2}^2 = \ip{K_z, K_z} = K_z(z) = \frac{1}{(1 - |z|^2)^2}.
\]
Thus
\begin{equation} \label{11}
k_z(w) = \frac{1 - |z|^2}{(1 - \bar{z}w)^2}
\end{equation}%
for $z, w \in D$.

Analytic automorphisms of the unit disk will play a key role
here.  For $z \in D$, let
$\phi_{z}$ be the M\"obius map on $D$ defined by
\begin{equation} \label{12}
\phi_{z}(w)=\frac{z-w}{1-\bar{z}w}.
\end{equation}%
Let $U_z \colon L^2_a \to L^2_a$ be the unitary
operator defined by
\begin{equation} \label{67}
U_z f = (f \circ \phi_z) {\phi_z}'.
\end{equation}%
To show that $U_z$ is indeed unitary, first make a change of variables in the
integral defining $\|U_zf\|_2$ to show that
$U_z$ is an isometry on $L^2_a$.  Next, a simple computation shows that ${U_z}^2$ is
the identity operator on $L^2_a$ (this holds because $\phi_z$ is its own inverse under
composition).  Being an invertible isometry, $U_z$ must be unitary.  Notice
that ${U_z}^* = {U_z}^{-1} = U_z$, so $U_z$ is actually a self-adjoint unitary
operator.

We will need two more simple properties of $U_z$.  First,
\begin{equation} \label{76}
U_z 1 = -k_z;
\end{equation}%
this follows from \eqref{11}, \eqref{12}, and \eqref{67}.  Second,
\begin{equation} \label{25}
U_z T_u U_z = T_{u \circ \phi_z}
\end{equation}%
for every $u \in L^\infty(D, dA)$; this is proved as Lemma~8 of \bibref{AxC}.  Thus
if
$u_1, \dots, u_n \in L^\infty(D, dA)$, then
\begin{equation} \label{26}
U_z T_{u_1} \dots T_{u_n} U_z =  T_{u_1 \circ \phi_z} \dots T_{u_n \circ \phi_z}
\end{equation}%
because we can write the operator on the left side as
\[
(U_z T_{u_1} U_z)(U_z T_{u_2} U_z) \dots (U_z T_{u_n} U_z)
\]
and then use \eqref{25}.

Next we compute the Berezin transform of a product of Toeplitz
operators.  The formula given by the following lemma will be used later when
we prove that $\tl{S} \in \U$ for every $S \in \T$ (Theorem~\ref{23}).

\begin{lemma} \label{17}
If\/ $u_1, \dots, u_n \in L^\infty(D, dA)$, then
\[
(T_{u_1} \dots T_{u_n})\til(z)
= \ip{ T_{u_1 \circ \phi_z} \dots T_{u_n \circ \phi_z} 1, 1}
\]
for every\/ $z \in D$.
\end{lemma}

\proof
Suppose $u_1, \dots, u_n \in L^\infty(D, dA)$ and $z \in D$.  
Then
\begin{align*}
(T_{u_1} \dots T_{u_n})\til(z)
&= \ip{ T_{u_1} \dots T_{u_n} k_z, k_z } \\
&= \ip{ T_{u_1} \dots T_{u_n} U_z 1, U_z 1} \\
&= \ip{ U_z T_{u_1} \dots T_{u_n} U_z 1, 1} \\
&= \ip{ T_{u_1 \circ \phi_z} \dots T_{u_n \circ \phi_z} 1, 1},
\end{align*}
where the first equality comes from the definition of the Berezin transform, the
second equality comes from
\eqref{76}, the third equality holds because $U_z$ is self-adjoint, and the last
equality comes from~\eqref{26}.
\qed

We will need to make extensive use of the \emph{maximal ideal space} of $H^{\infty}$,
which we denote by~$\M$.  We define $\M$ to be the set of  multiplicative linear maps
from $H^{\infty}$ onto the field of complex numbers.  With the weak-star topology,
$\M$ is a compact Hausdorff space. If $z$ is a point in  the unit disk~$D$, then point
evaluation at $z$ is a multiplicative linear functional on~$\M$.  Thus we can think of
$z$ as an element of
$\M$ and the unit disk $D$ as a subset of~$\M$.  Carleson's corona
theorem states that $D$ is dense in~$\M$.

Suppose $m \in \M$ and $z \mapsto \alpha_z$ is a mapping of $D$ into some
topological space~$E$.  Suppose also that $\beta \in E$.  The notation
\[
\lim_{z \to m} \alpha_z = \beta
\]
means (as you should expect) that for each open set $X$ in $E$ containing $\beta$,
there is an open set $Y$ in $\M$ containing $m$ such that $\alpha_z \in X$ for all
$z \in Y \cap D$.  Note that with this notation $z$ is always assumed to lie in $D$. 
We must deal with these nets rather than sequences because the topology of $\M$ is not
metrizable.

The Gelfand transform allows us to think of $H^{\infty}$ as contained in~$\CM$, the
algebra of continuous complex-valued functions on~$\M$.  By the Stone-Weierstrass
theorem, the set of finite sums of functions of the form $f\bar{g}$, with $f, g \in
H^\infty$, is dense in $\CM$, where $\CM$ is endowed with the usual supremum norm. 
Because $D$ is dense in $\M$, this supremum norm is the same as the usual supremum
norm over~$D$.  Thus we can identify $\CM$ with~$\U$, the closure in
$L^\infty(D, dA)$ of finite sums of functions of the form $f\bar{g}$, with $f, g \in
H^\infty$.

We will make frequent use of the identification discussed above of $\U$ with~$\CM$.  It
asserts that given a function $u \in \U$, which we normally think of as a function on
$D$, we can uniquely extend $u$ to a continuous complex-valued function on $\M$; this
extension to $\M$ is also denoted by~$u$.  Thus for $u \in \U$ and $m \in \M$, the
expression
$u(m)$ makes sense---it is the complex number defined by
\[
u(m) = \lim_{z \to m} u(z).
\]
Conversely, we will sometimes use the identification of $\U$ with $\CM$ to prove that
a function is in~$\U$.  Specifically, if $u$ is a continuous function on $D$ and we can
prove that $u$ extends to a continuous function on~$\M$, then we can conclude that
$u \in \U$.

For $m \in \M$, let $\phi_m \colon D \to \M$ denote the
\emph{Hoffman map}.  This is defined by setting
\[
\phi_m(w) = \lim_{z \to m} \phi_z(w)
\]
for $w \in D$; here we are taking a limit in~$\M$.  The existence of this limit,
as well as many other deep properties of $\phi_m$, was proved by Hoffman
\bibref{Hof}.  An exposition of Hoffman's results can also be found in \bibref{Gar},
Chapter~X.  We shall use, without further comment, Hoffman's result that
$\phi_m$ is a continuous mapping of $D$ into~$\M$.  Note that $\phi_m(0) = m$.

{\parfillskip=0pt plus1fil 
If $u \in \U$ and $m \in \M$, then $u \circ \phi_m$ makes sense as a continuous
function on~$D$, because $\phi_m$ maps $D$ into $\M$ and $u$ can be thought of as a
continuous function on $\M$, as we discussed above.  The next lemma provides the
crucial continuity that we will soon need.  Recall that a net of operators
$\{S_z\}_{z\in D} \subset \B(L^2_a)$ is said to converge to $S \in \B(L^2_a)$ in the
\emph{strong operator topology} as $z \to m$ if $\lim_{z \to m} S_z f = Sf$ for
every $f \in L^2_a$, where the last limit is taken in the norm in~$L^2_a$. 

}\begin{lemma} \label{18}
If \/ $u_1, \dots, u_n \in \U$, then
\begin{equation} \label{13}
\lim_{z \to m} T_{u_1 \circ \phi_z} \dots T_{u_n \circ \phi_z}
= T_{u_1 \circ \phi_m} \dots T_{u_n \circ \phi_m}
\end{equation}%
for every\/ $m \in \M$, where the limit is taken in the strong operator topology.
\end{lemma}

\proof
Fix $m \in \M$.  We will prove \eqref{13} by induction on $n$.  To
get the induction started, suppose
$n = 1$, so we consider a single function
$u \in \U$.  As $z \to m$, clearly $u \circ \phi_z$ converges to
$u \circ \phi_m$ pointwise on $D$.  Because the family of functions
\mbox{$\{u \circ \phi_z: z \in D\}$} is uniformly bounded, this
convergence is uniform on each compact subset of $D$ if $u$ happens to be analytic
on~$D$.  Thus the convergence is also uniform on each compact subset of $D$ if $u$
happens to be the product of an $H^\infty$ function and the complex conjugate
of an $H^\infty$ function.  Finite sums of such functions are dense in $\U$. 
Thus we can conclude that $u \circ \phi_z$ converges to
$u \circ \phi_m$ (as $z \to m$) uniformly on each compact subset of $D$ for
arbitrary $u \in \U$.  Fix $f \in L^2_a$.  Then
\begin{equation} \label{14}
\lim_{z \to m} \int_D | (u \circ \phi_z)(w) - (u \circ \phi_m)(w)|^2
|f(w)|^2\,dA(w) = 0,
\end{equation}%
because $D$, and hence the integral above, can be broken into two pieces---a
large compact subset of $D$ (on which $u \circ \phi_z$ converges uniformly to $u \circ
\phi_m$) and a set of small measure on which all the integrals
are small.  (We had to use uniform convergence on compact subsets of $D$ to prove
\eqref{14} because the Lebesgue dominated convergence theorem fails for nets, as
opposed to sequences.)  Because
\[
\lim_{z \to m} \| (u \circ \phi_z)f - (u \circ \phi_m)f \|_2 = 0,
\]
we have $\lim_{z \to m} \| T_{u \circ \phi_z} f - T_{u \circ \phi_m}f \|_2 =
0$, proving \eqref{13} in the case $n = 1$.

Now suppose that $u_1, \dots, u_n \in \U$ and that \eqref{13} holds when $n$ is
replaced by $n-1$.  For convenience, let
\[
S_z = T_{u_1 \circ \phi_z} \dots T_{u_{n-1} \circ \phi_z} \quad \text{and} \quad
S_m = T_{u_1 \circ \phi_m} \dots T_{u_{n-1} \circ \phi_m}.
\]
By our induction hypothesis, $\| S_z g - S_m g \|_2 \to 0$ as $z \to m$
for every $g \in L^2_a$.  Fix $f \in L^2_a$.  Then
\begin{align*}
\| T_{u_1 \circ \phi_z} \dots &T_{u_n \circ \phi_z} f
- T_{u_1 \circ \phi_m} \dots T_{u_n \circ \phi_m}f \|_2 \\
&= \| S_z T_{u_n \circ \phi_z} f - S_m T_{u_n \circ \phi_m} f \|_2 \\
&\le \| S_z \|_2 \|( T_{u_n \circ \phi_z} - T_{u_n \circ \phi_m}) f\|_2
+ \| (S_z - S_m)(T_{u_n \circ \phi_m} f) \|_2.
\end{align*}
Because $\|S_z\|_2$ is bounded by $\|u_1\|_\infty \dots \|u_n\|_\infty$, which is
independent of $z$, the first term in the last inequality above has limit $0$ as
$z \to m$ (by the $n=1$ case that we already proved).  Our induction hypothesis
implies that the second term in the last inequality above also has limit $0$ as
$z \to m$, completing the proof of~\eqref{13}.
\qed

Now we are ready to prove that the Berezin transform maps $\T$ into~$\U$.  
Most of the work needed to prove the theorem below was done in the last two
lemmas.  For the first time we will need to use the linearity of the Berezin transform
as well as its continuity: $\|\tl{S}\|_\infty \le \|S\|$ for all $S \in \B(L^2_a)$. 
We will also need to make use of the description of $\T$ as the closure in
$\B(L^2_a)$ of the set of finite sums of operators of the form
$T_{u_1} \dots T_{u_n}$, where $u_1, \dots, u_n \in \U$.
 
\begin{theorem} \label{23}
If\/ $S \in \T$, then\/ $\tl{S} \in \U$.  Furthermore,
if\/ $u_1, \dots, u_n \in \U$, then
\begin{equation} \label{15}
(T_{u_1} \dots T_{u_n})\til(m)
= \ip{ T_{u_1 \circ \phi_m} \dots T_{u_n \circ \phi_m} 1, 1}
\end{equation}%
for every\/ $m \in \M$.
\end{theorem}

\proof
If $u_1, \dots, u_n \in \U$, then Lemma~\ref{17} and Lemma~\ref{18} show that
$(T_{u_1} \dots T_{u_n})\til$ extends to be a continuous function on $\M$ and that the
extension is given by \eqref{15}.  Thus the Berezin transform maps sums of operators
of the form $T_{u_1} \dots T_{u_n}$, where each $u_j \in \U$, into $\U$.  The
linearity and continuity of the Berezin transform now imply that the Berezin
transform also maps $\T$ into~$\U$.
\qed

A multiplicative linear function $m \in \M$ is called a
\emph{one-point part} if $\phi_m$ is a constant map.  In other words, $m$ is a
one-point part if $\phi_m(w) = m$ for every $w \in D$.  The set of all one-point
parts is denoted by $\M_{1}$.  As is well known, $\M_{1}$ is a
closed subset of $\M$ that properly contains the Shilov boundary of
$H^\infty$ (in particular, $\M_{1}$ is not the empty set).  Actually
$\M_{1}$ should be thought of as a small subset of $\M
\setminus D$, as the complement of $\M_{1}$ in
$\M \setminus D$ is dense in $\M \setminus D$.

The following corollary shows how to compute the Berezin transform on $\M_1$ of a
finite product of Toeplitz operators with symbols in~$\U$.

\begin{corollary} \label{19}
If\/ $u_1, \dots, u_n \in \U$, then
\begin{equation} \label{16}
(T_{u_1} \dots T_{u_n})\til(m)
= u_1(m) \dots u_n(m)
\end{equation}%
for every\/ $m \in \M_1$.
\end{corollary}

\proof
Suppose $u_1, \dots, u_n \in \U$ and $m \in \M_1$.  Because $m \in \M_1$, each
function $u_j \circ \phi_m$ is a constant function equal to the constant $u_j(m)$. 
Thus each of the Toeplitz operators in \eqref{15} has constant symbol, reducing
\eqref{15} to the desired equation \eqref{16}.
\qed

{\parfillskip=0pt plus1fil
The next corollary shows that the Berezin transform is multiplicative on~$\M_1$.

}
\begin{corollary} \label{20}
If\/ $R, S \in \T$, then
\[
(RS)\til(m) = \tl{R}(m) \tl{S}(m)
\]
for every\/ $m \in \M_1$.
\end{corollary}

\proof
If $R, S$ are each products of Toeplitz operators with symbols in~$\U$, then
the desired result follows from Corollary~\ref{19}.  The proof is completed by
recalling that sums of such operators are dense in~$\T$.
\qed

Recall that the \emph{commutator ideal} $\CT$ is the smallest closed, two-sided ideal
of $\T$ containing all operators of the form $RS - SR$, where $R, S \in \T$.
In a remarkable theorem, McDonald and Sundberg (\bibref{McS}, Theorem~6; also see
\bibref{Sun} for another proof) showed that
$\T/\CT$ is isomorphic, as a $C^*$-algebra, to~$C(\M_1)$.  More precisely, they showed
that the map
$u \mapsto T_u + \CT$  is a surjective homomorphism of $\U$ onto $\T/\CT$, with kernel
$\{ u\in \U : u|_{\M_1} =0 \}$.  Rephrased again, the McDonald-Sundberg theorem
states that each $S \in \T$ can be written in the form $S = T_u + R$ for some
$u \in \U$ and some $R \in \CT$.  Furthermore, if $u \in \U$, then $T_u \in \CT$ if
and only if $u|_{\M_1} = 0$.  These results account for the importance of
understanding the commutator ideal~$\CT$.  We now describe
$\CT$ in terms of Berezin transforms.

\begin{theorem} \label{22}
Suppose\/ $S \in \T$.  Then\/ $S$ is in the commutator ideal\/ $\CT$ if and only if\/
$\tl{S}|_{\M_1} = 0$.
\end{theorem}

\proof
The commutator ideal $\CT$ is the norm closure of the set of finite sums of operators
of the form $S_1(S_2 S_3 - S_3 S_2)S_4$, where $S_1, S_2, S_3, S_4 \in \T$.  By
Corollary~\ref{20}, each such operator has a Berezin transform that vanishes
on~$\M_1$.  Thus if $S \in \CT$, then $\tl{S}|_{\M_1} = 0$, proving one direction
of the theorem.

To prove the other direction, suppose $\tl{S}|_{\M_1} = 0$. By the McDonald-Sundberg
theorem, we can  write  $S =T_u + R$ for some $u \in \U$ and $R\in \CT$. Thus
\begin{align*}
0 &= \tl{S}|_{\M_1} \\
&= \wt{T_u}|_{\M_1} + \tl{R}|_{\M_1} \\
&= u|_{\M_1} + \tl{R}|_{\M_1} \\
&= u|_{\M_1},
\end{align*}
where the third equality comes from Corollary~\ref{19} and the fourth equality holds by
the direction of this theorem that we have already proved.  The McDonald-Sundberg
theorem now tells us that $T_u \in \CT$ (because $u|_{\M_1} = 0$). Therefore
$S \in \CT$, completing the proof.
\qed

Given an operator $S \in \T$, the McDonald-Sundberg theorem tells us that $S$ can be
written in the form $S = T_u + R$ for some $u \in \U$ and $R\in \CT$.  The choice of
$u$ is not unique, as it can be perturbed by any function in $\U$ that vanishes
on~$\M_1$.  However, we now show that there is a canonical choice of~$u$, namely the
Berezin transform of~$S$.  The corollary below states that the decomposition
\[
S = T_{\tl{S}} + (S - T_{\tl{S}})
\]
satisfies the requirements of the McDonald-Sundberg theorem, because the term in
parentheses is in~$\CT$.

\begin{corollary} \label{66}
If\/ $S \in \T$, then\/ $S - T_{\tl{S}} \in \CT$.
\end{corollary}

\proof
Suppose $S \in \T$.  Then by Theorem~\ref{23}, $\tl{S} \in \U$.  If $m \in \M_1$, then
using Corollary~\ref{19} (with $n=1$) we get
\[
(S - T_{\tl{S}})\til(m) = \tl{S}(m) - \tl{S}(m) = 0.
\]
In other words, $(S - T_{\tl{S}})\til|_{\M_1} = 0$.  Thus Theorem~\ref{22} implies
that $S - T_{\tl{S}} \in \CT$, completing the proof.
\qed

The next lemma will allow us to translate results about $\M_1$, a rather abstract
object, into results about nontangential behavior on the unit disk~$D$.  When we
refer to ``almost every point of $\partial D$'', we mean with respect to the usual
linear Lebesgue (arc length) measure on~$\partial D$.

\begin{lemma} \label{21}
If\/ $u \in \U$ and\/ $u|_{\M_1} = 0$, then\/ $u$ has nontangential limit\/ $0$ at
almost every point of\/ $\partial D$.
\end{lemma}

\proof
As is well known, every function in $H^\infty$ has a nontangential limit at almost
every point of~$\partial D$.  Thus every finite sum of functions of the form
$f\bar{g}$, where $f, g \in H^\infty$, has a nontangential limit at almost
every point of~$\partial D$.  Hence any function on $D$ that is the uniform limit of a
sequence of such functions also has a nontangential limit at almost
every point of~$\partial D$ (this holds because the union of a countable
collection of sets of measure $0$ has measure $0$).  In other words, every function in
$\U$ has a nontangential limit at almost every point of~$\partial D$.

Suppose $u \in \U$ and $u|_{\M_1} = 0$.  Define a function $u^*$ (almost
everywhere) on $\partial D$ by letting $u^*(\lambda)$ equal the nontangential limit
of $u$ at $\lambda \in \partial D$.  Let $X \subset \M$ denote the Shilov boundary of
$H^\infty$.  By Theorem~11 of Axler and Shields's paper \bibref{ASh}, the essential
range of $u^*$ on $\partial D$ equals~$u(X)$.  However, $X$ is contained in $\M_1$, so
we conclude that the essential range of $u^*$ on $\partial D$ is just $\{0\}$.  Thus
$u^*$ equals $0$ almost everywhere on $\partial D$.  Hence $u$ has nontangential
limit $0$ at almost every point of $\partial D$.
\qed

{\parfillskip=0pt plus1fil
Now we can prove that the Berezin transform of each operator
in the commutator ideal of $\T$ has nontangential limit $0$ almost everywhere
on~$\partial D$.

}\begin{corollary} \label{77}
If\/ $S \in \CT$, then\/ $\tl{S}$ has nontangential limit\/ $0$ at
almost every point of\/ $\partial D$.
\end{corollary}

\proof
Combine Theorem~\ref{22} and Lemma~\ref{21} to obtain the desired result.
\qed

The converse of the corollary above is false.  To see this, let $u$ be a function in
$\CM$ that equals $0$ on the Shilov boundary of $H^\infty$ but that is not
identically $0$ on~$\M_1$.  Then $\wt{T_u}$ equals $0$ on the Shilov boundary
of~$H^\infty$ (by Corollary~\ref{19}).  The proof of Lemma~\ref{21} thus shows that
$\wt{T_u}$ has nontangential limit $0$ at almost every point of~$\partial D$. 
However, $\wt{T_u}$ is not identically~$0$ on $\M_1$ (by Corollary~\ref{19}) and thus
$T_u$ is not in $\CT$ (by Theorem~\ref{22}), providing the desired example.

The next lemma will be used in the proof of Proposition~\ref{31}.  In equation
\eqref{30} below,
$1$ denotes the constant function on
$D$ (the function that maps $z$ to $1$) and $z$ denotes the identity function on $D$
(the function that maps $z$ to $z$).

\begin{lemma} \label{29}
If\/ $S \in \B(L^2_a)$, then\/ $\tl{S}$ is real analytic on\/ $D$ and
\begin{equation} \label{30}
(\Delta \tl{S})(0) = 16 \ip{Sz,z} - 8 \ip{S1,1}.
\end{equation}%
\end{lemma}

\proof
Let $S \in \B(L^2_a)$.  Define a complex-valued function $F$ on $D \times D$ by
\[
F(w, z) = \ip{SK_{\bar{w}}, K_z}
\]
for $w, z \in D$.  Note that here we are using the unnormalized reproducing kernels. 
For fixed $w \in D$, the function $SK_{\bar{w}}$ is in $L^2_a$, and hence is analytic
on~$D$.  Because $F(w, z) = (SK_{\bar{w}})(z)$, this implies that $F(w, z)$ is
analytic in $z$ for fixed~$w$.  Similarly, for fixed $z \in D$, the function $S^*K_z$
is in $L^2_a$, and hence is analytic on~$D$.  Because
$F(w, z) = \overline{(S^*K_z)(\bar{w})}$, this implies that $F(w, z)$ is analytic in
$w$ for fixed~$z$.  Because $F$ is analytic in each variable separately, we conclude
that $F$ is holomorphic on $D \times D$.
Clearly $\tl{S}(z) = (1 - |z|^2)^2 F(\bar{z}, z)$.  Because $F$ is holomorphic on
$D \times D$, this implies that $\tl{S}$ is real analytic on~$D$, as desired.

To prove \eqref{30}, we first express the explicit formula \eqref{11} for the normalized
reproducing kernel as a power series:
\[
k_z(w) = (1 - |z|^2) \sum_{j=0}^\infty (j+1) \bar{z}^j w^j.
\]
Thus
\begin{align}
{}\quad\tl{S}(z) &= \ip{ Sk_z, k_z } \notag \\
&= (1 - |z|^2)^2
\sum_{j,n=0}^\infty (j+1) (n+1) \ip{ Sw^j, w^n } \bar{z}^j z^n \label{79} \\
&= (1 - 2 z \bar{z} + z^2 \bar{z}^2)
\sum_{j, n=0}^\infty (j+1) (n+1) \ip{ Sw^j, w^n } \bar{z}^j z^n \label{48} \\
&= \sum_{j, n=0}^\infty a_{j,n} \bar{z}^j z^n, \label{47}
\end{align}
where the coefficients $a_{j,n}$ could be computed explicitly.  Note that
\begin{align*}
(\Delta \tl{S})(0) &= 4 \frac{\partial^2 \tl{S}}{\partial \bar{z}\,\partial
z}(0)
\\ &= 4 a_{1,1},
\end{align*}
where the last equation follows from \eqref{47}.  From \eqref{48} we see that
\[
a_{1,1} = 4 \ip{Sw, w} - 2 \ip{ S1, 1}.
\]
The proof of \eqref{30} is completed by combining the last two equations and replacing
the independent variable $w$ above (denoting the identity function) with the more
common symbol $z$.
\qed

We note for later use that \eqref{79} implies that an operator
$S \in \B(L^2_a)$ is uniquely determined by its Berezin transform.  To see this,
suppose $S \in \B(L^2_a)$ and $\tl{S}$ is identically $0$ on~$D$.  We
need to show that $S=0$.  Differentiating the infinite sum in \eqref{79} $n$ times with
respect to $z$ and
$j$ times with respect to $\bar{z}$ and then evaluating at $z = 0$ shows that
\mbox{$\ip{Sw^j, w^n} = 0$} for all nonnegative integers $j$ and $n$.  Because finite
linear combinations of
$\{w^j: j \ge 0\}$ are dense in~$L^2_a$, this implies that $S=0$, as desired.

For $S \in \B(L^2_a)$ and $z \in D$, define $S_z \in \B(L^2_a)$ by
\[
S_z = U_z S U_z.
\]
Recall that $U_z$ was defined by equation \eqref{67}.  Note that if $S$ is a finite
product of Toeplitz operators, then a formula for
$S_z$ is given by~\eqref{26}.  The next lemma shows us how to define $S_m$ for each
$m \in \M$ in a manner consistent with the definition just given for $S_z$.  This
operator $S_m$ plays an important role in Proposition~\ref{31}, where it is used in the
proofs of parts (b) and (c) even though it is not explicitly mentioned in the
statements of those results.

\begin{lemma} \label{28}
If\/ $S \in \T$ and\/ $m \in \M$, then there exists\/ $S_m \in \T$ such that
\[
\lim_{z \to m} S_z = S_m,
\]
where the limit is taken in the strong operator topology.  If\/
$S = T_{u_1} \dots T_{u_n}$, where\/ $u_1, \dots, u_n \in \U$, then\/
$S_m = T_{u_1 \circ \phi_m} \dots T_{u_n \circ \phi_m}$.
\end{lemma}

\proof
Fix $m \in \M$.  First suppose $S = T_{u_1} \dots T_{u_n}$, where
$u_1, \dots, u_n \in \U$.  Then
$S_z = T_{u_1 \circ \phi_z} \dots T_{u_n \circ \phi_z}$ for every $z \in D$, as
we saw in~\eqref{26}.  Thus, by Lemma~\ref{18},
\[
\lim_{z \to m} S_z = T_{u_1 \circ \phi_m} \dots T_{u_n \circ \phi_m}.
\]
To prove that the operator on the right side of this equation is in $\T$, we must show
that $u \circ \phi_m \in \U$ whenever $u \in \U$.  Clearly this holds if
$u \in H^\infty$, because then $u \circ \phi_m \in H^\infty$.  Taking complex
conjugates and then products, we have that $u \circ \phi_m \in \U$ for all $u$ of the
form $f\bar{g}$, where $f, g \in H^\infty$.  Finite sums of such functions are dense in
$\U$, showing that $u \circ \phi_m \in \U$ for all $u \in \U$, as desired.  This
completes the proof of the lemma when $S$ has the special form $T_{u_1} \dots
T_{u_n}$, where $u_1, \dots, u_n \in \U$.

Now suppose $S \in \T$.  Fix $f \in L^2_a$.  We must prove that
$\lim_{z \to m} S_z f$ exists in~$L^2_a$.  To do this, suppose
$\epsilon > 0$.  Then there is an operator $R$ that is a finite sum of operators of the
form $T_{u_1} \dots T_{u_n}$, where each $u_j \in \U$, such that
$\|S - R\| \le \epsilon$.  Thus $\|S_z f - R_z f\|_2 \le \epsilon
\|f\|_2$.  {}From the paragraph above, we know that $R_z f$ converges (as $z
\to m$) to a function $R_m f$.  Thus
\[
\limsup_{z \to m} \| S_z f - R_m f\|_2 \le \epsilon \|f\|_2.
\]
Thus
\[
\limsup_{z, w \to m} \| S_z f - S_w f\|_2 \le 2\epsilon \|f\|_2.
\]
Because $\epsilon$ is an arbitrary positive number, this means that $S_z f$ is a
Cauchy net in $L^2_a$ (as $z \to m$).  However, $L^2_a$ is complete, and so
this Cauchy net must converge, as desired.

{}From the first paragraph of this proof, we know that $S_m \in \T$ for all $S$
in a dense subset of $\T$.  The mapping $S \mapsto S_m$ is continuous (in the operator
norm), so $S_m$ must be in $\T$ for all $S \in \T$, completing the proof.
\qed

A function $u \in \U$ is said to be \emph{real analytic} on $\M$ if $u \circ \phi_m$
is real analytic on $D$ for every $m \in \M$.  The next proposition tells us that the
Berezin transform of any operator $S \in \T$ is real analytic on~$\M$.  Furthermore,
for $m \in \M$ we get a formula for computing the Laplacian at $0$ of
$\tl{S} \circ \phi_m$.  These results will be used in the next section of this
paper.

\begin{proposition} \label{31}
Suppose\/ $S \in \T$.  Then
\begin{subtheorem}
\item
$\tl{S} \circ \phi_m = \wt{S_m}$ for every\/ $m \in \M$;

\item
$\tl{S}$ is real analytic on\/ $\M$;

\item
$\displaystyle \bigl(\Delta (\tl{S} \circ \phi_m)\bigr)(w) =
\lim_{z \to m} \frac{(1 - |\phi_z(w)|^2)^2 (\Delta \tl{S})(\phi_z(w))}{(1 -
|w|^2)^2}$\\
for every\/ $w \in D, m \in \M$.
\end{subtheorem}
\end{proposition}

\proof
We begin by deriving a useful formula.  Suppose $w, z \in D$.  If $f \in L^2_a$, then
\begin{align*}
\ip{f, U_z K_w} &= \ip{U_z f, K_w} \\
&= (U_z f)(w) \\
&= (f \circ \phi_z)(w){\phi_z}'(w) \\
&= \ip{f, \overline{{\phi_z}'(w)} K_{\phi_z(w)}}.
\end{align*}
Thus $U_z K_w = \overline{{\phi_z}'(w)} K_{\phi_z(w)}$.  Rewriting this in terms of
the normalized reproducing kernels, we have
\begin{equation} \label{27}
U_z k_w = \alpha k_{\phi_z(w)}
\end{equation}%
for some complex constant $\alpha$.  Without doing a computation, we know that
$|\alpha| = 1$, because
\mbox{$\|k_w\|_2 = \|k_{\phi_z(w)}\|_2 = 1$} and $U_z$ is unitary.

For the rest of the proof, fix $m \in \M$.  To prove (a), fix $w \in D$.  If
$z \in D$, then
\begin{align*}
\tl{S}\bigl(\phi_z(w)\bigr) &= \ip{Sk_{\phi_z(w)}, k_{\phi_z(w)}} \\
&= \ip{SU_z k_w, U_z k_w} \\
&= \ip{U_z S U_z k_w, k_w},
\end{align*}
where the second equality comes from \eqref{27} along with the extra information that
$|\alpha| = 1$.  Taking limits of the first and last terms above as $z \to m$,
we get $\tl{S}\bigl(\phi_m(w)\bigr) = \ip{S_m k_w, k_w}$.  Thus (a) holds.

To prove (b), recall that $S_m \in \T$ (by Lemma~\ref{28}).  Thus $\wt{S_m}$ is real
analytic on $D$ (by Lemma~\ref{29}).  Now (a) shows that $\tl{S} \circ \phi_m$ is
real analytic on $D$.  We thus conclude that $\tl{S}$ is real analytic on $\M$,
completing the proof of~(b).

To prove (c), fix $w \in D$.  Then
\begin{align*}
(\Delta \tl{S})(w)
&= \frac{\bigl(\Delta (\tl{S} \circ \phi_w)\bigr)(0)}{(1-|w|^2)^2} \\[9pt]
&= \frac{(\Delta \wt{S_w})(0)}{(1-|w|^2)^2} \\[9pt]
&= \frac{16\ip{SU_wz, U_wz} - 8\ip{SU_w1, U_w1}}{(1-|w|^2)^2},
\end{align*}
where the first equality follows from a
standard calculation, the second equality comes from~(a), and the third inequality
comes from Lemma~\ref{29}.  The equation above shows that the map
$S \mapsto (\Delta \tl{S})(w)$ is a continuous linear functional on $\B(L^2_a)$ with
respect to the strong operator topology on~$\B(L^2_a)$.  Now fix $m \in \M$.  Then
\begin{align*}
\bigl(\Delta (\tl{S} \circ \phi_m)\bigr)(w)
&= (\Delta \wt{S_m})(w) \\
&= \lim_{z \to m} (\Delta \wt{S_z})(w) \\
&= \lim_{z \to m} \bigl(\Delta (\tl{S} \circ \phi_z)\bigr)(w) \\
&= \lim_{z \to m} \frac{(1 - |\phi_z(w)|^2)^2 (\Delta \tl{S})(\phi_z(w))}{(1 -
|w|^2)^2},
\end{align*}
where the first equality holds by (a), the second equality holds by the
continuity discussed earlier in this paragraph and Lemma~\ref{28}, the
third inequality holds by (a), and the fourth inequality holds by a standard
calculation.  This completes the proof of~(c).
\qed

\section{Boundary Behavior of the Berezin Transform on $\U$}

In this section we will study the boundary behavior of the Berezin transform on
elements of~$\U$.  Recall that if $u \in L^\infty(D, dA)$, then the Berezin transform
$\tl{u}$ is the function on $D$ defined by $\tl{u} = \wt{T_u}$.  This definition
leads to the explicit formula~\eqref{49}.

Our next result states that if $u$ is in $\U$, then so is $\tl{u}$.  If $u \in \U$,
then $u \circ \phi_{m}$ is bounded and continuous on~$D$, so the integral appearing in
the proposition below makes sense.

\begin{proposition} \label{9}
The Berezin transform maps\/ $\U$ into\/ $\U$.  Furthermore, if\/ $u \in \U$, then
\begin{equation} \label{7}
\tl{u}(m)=\int_{D} (u \circ \phi_{m})(w)\,dA(w)
\end{equation}%
for every\/ $m \in M$.
\end{proposition}

\proof
Suppose $u \in \U$.  By definition, $\tl{u} = \wt{T_u}$.  Theorem~\ref{23} thus
tells us that $\tl{u} \in \U$.  Furthermore, from \eqref{15}, which is used in the
second equality below, we have
\begin{align*}
\tl{u}(m) &= \wt{T_u}(m) \\
&= \ip{T_{u \circ \phi_m} 1, 1} \\
&= \ip{ P(u \circ \phi_m), 1} \\
&= \ip{ u \circ \phi_m, 1} \\
&=\int_{D} (u \circ \phi_{m})(w)\,dA(w),
\end{align*}
for every $m \in M$, as desired.
\qed

If $u$ is a bounded harmonic function on $D$, then so is $u \circ \phi_z$ for each
$z \in D$.  The mean value property and \eqref{7} then imply that
$\tl{u}(z) = (u \circ \phi_z)(0) = u(z)$ for each $z \in D$.  In other words, every
harmonic function equals its Berezin transform.

The next corollary shows that a function in $\U$ and its Berezin transform agree on
the set of one-point parts.

\begin{corollary} \label{24}
If\/ $u \in \U$, then\/ $\tl{u}|_{\M_1} = u|_{\M_1}$.
\end{corollary}

\proof
Suppose $u \in \U$ and $m \in \M_1$.  Then $(u \circ \phi_m)(w) = u(m)$
for every $w \in D$ (recall that $m \in \M_1$ implies that $\phi_m$ is a constant map
on~$D$).  Thus \eqref{7} shows that $\tl{u}(m) = u(m)$, as desired.
\qed

The next corollary shows that a function in $\U$ and its Berezin transform have the
same nontangential limits almost everywhere on $\partial D$ (recall from the proof of
Lemma~\ref{21} that every function in $\U$ has nontangential limits almost everywhere
on~$\partial D$).

\begin{corollary} \label{42}
If\/ $u \in \U$, then\/ $\tl{u} - u$ has nontangential limit\/ $0$ at almost every
point of\/ $\partial D$.
\end{corollary}

\proof
Suppose $u \in \U$.  Then $\tl{u} - u \in \U$ (from Proposition~\ref{9}) and
$(\tl{u} - u)|_{\M_1} = 0$ (from Corollary~\ref{24}).  Lemma~\ref{21} now gives the
desired result.
\qed

For $u \in L^\infty(D, dA)$, the \emph{Hankel operator} with symbol $u$ is the
operator $H_u$ from $L^2_a$ to $L^2(D,dA) \ominus L^2_a$ defined by
$H_u f = (1 - P)(uf)$.  The next corollary shows that Toeplitz and Hankel operators
with symbol in $\U$ behave nicely on normalized reproducing kernels corresponding to
a net of points converging to a one-point part.
   
\begin{corollary} \label{10}
If\/ $u \in \U$, then
\[
\lim_{z \to m}\|T_{u-u(m)}k_z\|_{2}=0
\]
and
\[
\lim_{z \to m}\|H_{u}k_z\|_{2}=0
\]
for every\/ $m \in \M_{1}$.
\end{corollary}

\proof
Suppose $u \in \U$ and $m \in \M_{1}$.  We claim that
\begin{equation} \label{8}
\lim_{z \to m} \|(u - u(m))k_z\|_2 = 0.
\end{equation}%
Once this is proved, the proof will be done, because
\[
\|T_{u-u(m)}k_z\|_{2} = \|P\bigl((u - u(m))k_z\bigr)\|_2
\le \|(u - u(m))k_z\|_2
\]
and
\[
\|H_u k_z\|_{2} = \|H_{u-u(m)}k_z\|_{2}
= \|(1-P)\bigl((u - u(m))k_z\bigr)\|_2
\le \|(u - u(m))k_z\|_2.
\]

To prove \eqref{8}, note that
\begin{align*}
\lim_{z \to m} {\|(u - u(m))k_z\|_2}^2
&= \lim_{z \to m} \int_D |u(w) - u(m)|^2
|k_z(w)|^2\,dA(w) \\
&= \lim_{z \to m} (|u - u(m)|^2)\til(z) \\
&= (|u - u(m)|^2)\til(m),
\end{align*}
where the last equality holds because $(|u - u(m)|^2)\til$ can be thought of as a
continuous function on $\M$ (by Proposition~\ref{9}).  The function
$|u - u(m)|^2$ is in $\U$, so it equals its Berezin transform on
$\M_{1}$ (by Corollary~\ref{24}).  In particular, because $m \in \M_{1}$ and $|u -
u(m)|^2$ equals $0$ at $m$, the last quantity above equals $0$, completing the proof.
\qed

In the next corollary we once again translate a statement involving $\M_1$ into a
more concrete statement.

\begin{corollary} \label{43}
If\/ $u \in \U$, then the functions
\[
z \mapsto \|T_{u-u(z)}k_z\|_{2} \quad \text{and} \quad z \mapsto \|H_{u}k_z\|_{2}
\]
have nontangential limits\/ $0$ at almost every point of\/
$\partial D$.
\end{corollary}

\proof
Suppose $u \in \U$.  As in the proof of the previous corollary, we need only show that
$\|(u - u(z))k_z\|_2$ has nontangential limit $0$ at almost every point
of~$\partial D$.  To do this, note that
\begin{align*}
{\|(u - u(z))k_z\|_2}^2 &= \int_D |u(w) - u(z)|^2 |k_z(w)|^2\,dA(w) \\
&= \int_D \bigl(|u(w)|^2 - 2 \re(\overline{u(z)}u(w)) + |u(z)|^2\bigr)
|k_z(w)|^2\,dA(w) \\
&= \wt{|u|^2}(z) - 2 \re (\overline{u(z)}\tl{u}(z)) + |u(z)|^2
\end{align*}
for $z \in D$.  The equation above, along with Proposition~\ref{9}, shows that the
function $z \mapsto {\|(u - u(z))k_z\|_2}^2$ is in $\U$.  The equation above, along
with Corollary~\ref{24}, shows that the function $z \mapsto {\|(u - u(z))k_z\|_2}^2$
is $0$ on~$\M_1$.  Lemma~\ref{21} now gives the desired result.
\qed

The next two lemmas will be useful in proving Theorem~\ref{36}, which is the main
result of this section.  The formula for
$(\Delta \tl{u})(0)$ given by the first lemma below could be proved by differentiating
twice under the integral in the explicit formula for
$\tl{u}$ obtained from \eqref{49} and \eqref{11}.  However we have avoided that
computation in our proof by using the formula given by Lemma~\ref{29}.

\begin{lemma} \label{33}
If\/ $u \in L^\infty(D, dA)$, then\/ $\tl{u}$ is real analytic
on\/
$D$ and
\[
(\Delta \tl{u})(0) = 8 \int_D u(z)(2|z|^2 - 1)\,dA(z).
\]
\end{lemma}

\proof
Suppose $u \in L^\infty(D, dA)$.  Then $\tl{u} = \wt{T_u}$, and hence Lemma~\ref{29}
implies that $\tl{u}$ is real analytic on~$D$.  {}From Lemma~\ref{29} we also have
\begin{align*}
(\Delta \tl{u})(0) &= (\Delta \wt{T_u})(0) \\
&= 16 \ip{T_u z, z} - 8 \ip{T_u 1, 1} \\
&= 16 \ip{u z, z} - 8 \ip{u, 1} \\
&= 8 \int_D u(z)(2|z|^2 - 1)\,dA(z),
\end{align*}
completing the proof.
\qed

The next lemma provides information about the Berezin transforms of functions
analogous to the information about the Berezin transforms of operators provided by
Proposition~\ref{31}.

\begin{lemma} \label{32}
Suppose\/ $u \in \U$.  Then
\begin{subtheorem}
\item
$\tl{u} \circ \phi_m = (u \circ \phi_m)\til$ for every\/ $m \in \M$;

\item
$\tl{u}$ is real analytic on\/ $\M$;

\item
$\displaystyle \bigl(\Delta (\tl{u} \circ \phi_m)\bigr)(w) =
\lim_{z \to m} \frac{(1 - |\phi_z(w)|^2)^2 (\Delta \tl{u})(\phi_z(w))}{(1 -
|w|^2)^2}$ \\
for every\/ $w \in D, m \in \M$.
\end{subtheorem}
\end{lemma}

\proof
Suppose $m \in \M$.  Then
\begin{align*}
\tl{u} \circ \phi_m &= \wt{T_u}  \circ \phi_m \\
&= \bigl((T_u)_m\bigr)\til \\
&= \bigl(T_{u \circ \phi_m}\bigr)\til \\
&= (u \circ \phi_m)\til,
\end{align*}
where the second equality comes from Proposition~\ref{31}(a) and the third equality
comes from the second statement in Lemma~\ref{28}.  The equation above shows that (a)
holds.

Because $\tl{u} = \wt{T_u}$, (b) and (c) follow immediately from parts (b) and (c) of
Proposition~\ref{31}.
\qed

Now we turn to the question of describing the functions $u \in \U$ such that
$\lim_{z \to \partial D} \tl{u}(z) - u(z) = 0$.  Because the disk $D$ is dense in
$\M$, this is easily seen to be equivalent to the question of describing the
functions $u \in \U$ such that $\tl{u}$ equals $u$ on $\M \setminus D$.  We have
seen that $\tl{u}$ equals $u$ on~$\M_1$ for every $u \in \U$ (Corollary~\ref{24});
now we are asking when equality holds on the larger set $\M \setminus D$.
As motivation for our answer, recall that we pointed out earlier that every bounded
harmonic function equals its Berezin transform.  The converse also holds, so a function in $L^\infty(D, dA)$ equals its Berezin transform
if and only if it is harmonic (for proofs of this deep result, see the papers by
Engli\v s \bibref{Eng} or Ahern, Flores, and Rudin \bibref{AFR}).  Thus we might guess
that a function $u \in \U$ equals $\tl{u}$ on $\M \setminus D$ if and only if $u$ is
harmonic on $\M \setminus D$ (whatever that means).  As we will see, this turns out to
be correct if we define the notion of harmonic on $\M \setminus D$ in terms of the
parameterizations given by the Hoffman maps.

Motivated by the paragraph above, we define $\hop$ (which stands for ``harmonic on
parts'') to be the set of functions $u \in \U$ such that $u \circ \phi_m$ is harmonic
on $D$ for every $m \in \M \setminus D$.  Every bounded harmonic function on $D$ is in
$\hop$.  (Proof: If $u$ is a bounded harmonic function on $D$, then so is
$u \circ \phi_z$ for every $z \in D$.  Now $u \circ \phi_m(w) = \lim_{z \to m}u \circ
\phi_z(w)$ for every $w \in D, m \in \M$.  Because the pointwise limit of any uniformly
bounded net of harmonic functions is harmonic, we conclude that $u \circ \phi_m$ is
harmonic, as desired.)  Every function in $C(\bar{D})$ is also in $\hop$ (because if
$u \in C(\bar{D})$ and $m \in \M \setminus D$, then $u \circ \phi_m$ is a constant
function on~$D$).

The next theorem gives several conditions on a function $u \in \U$ that are
equivalent to having $\lim_{z \to \partial D} \tl{u}(z) - u(z) = 0$.  Note
that condition~(h) in the theorem below would not make sense for an arbitrary
$u \in \U$ (because functions in $\U$ need not even be differentiable on~$D$).  Even
for a function $u \in \U$ that is differentiable on $D$, there is no obvious connection
between the derivatives of $u$ on $D$ and derivatives of the functions
$u \circ \phi_m$ for $m \in \M \setminus D$.  This helps explain the extra hypothesis
required below for the applicability of condition~(h).

\begin{theorem} \label{36}
Suppose\/ $u \in \U$.  Then the following are equivalent:
\begin{subtheorem}
\item
${\displaystyle\lim_{z \to \partial D} \tl{u}(z) - u(z) = 0}$;

\item
$\tl{u} = u$ on\/ $\M \setminus D$;

\item
$u \in \hop$;

\item
$\tl{u} \in \hop$;

\item
$T_{\tl{u}} - T_u$ is a compact operator;

\item
${\displaystyle\lim_{z \to \partial D}
\int_D (u \circ \phi_z)(w) (2|w|^2 - 1)\,dA(w) = 0}$;

\item
${\displaystyle\lim_{z \to \partial D} (1 - |z|^2)^2 (\Delta \tl{u})(z) = 0}$.
\end{subtheorem}
If\/ $u$ is a finite sum of functions of the form\/ $u_1 \dots u_n$, where each\/ $u_j$
is a bounded harmonic function on\/ $D$, then the conditions above are also equivalent
to the condition below:
\begin{subtheorem} \setcounter{subtheoremc}{7}
\item
${\displaystyle\lim_{z \to \partial D} (1 - |z|^2)^2 (\Delta u)(z) = 0}$.
\end{subtheorem}
\end{theorem}

\proof
As is well known, the equivalence of (a) and (b) follows from the corona theorem.

Suppose (b) holds, so $\tl{u} = u$ on $\M \setminus D$.  Let
$m \in \M \setminus D$.  Then
\[
(u \circ \phi_m)\til = \tl{u} \circ \phi_m = u \circ \phi_m,
\]
where the first equality comes from Lemma~\ref{32}(a) and the second equality comes
from our hypothesis~(b).  The equation above says that $u \circ \phi_m$ is a function
in $L^\infty(D,dA)$ that equals its Berezin transform.  As we discussed earlier,
Engli\v s \bibref{Eng} and Ahern, Flores, and Rudin \bibref{AFR} proved that only
harmonic functions equal their Berezin transforms.  Thus $u \circ \phi_m$ is harmonic. 
Because $m$ was an arbitrary element of $\M \setminus D$, this implies that $u \in
\hop$.  Thus (b) implies (c).

Now suppose (c) holds, so $u \in \hop$.  If $m \in \M \setminus D$, then
\[
\tl{u} \circ \phi_m = (u \circ \phi_m)\til = u \circ \phi_m,
\]
where the first equality comes from Lemma~\ref{32}(a) and the second equality
holds because $u \circ \phi_m$ is harmonic.  The equation above shows that
$\tl{u} \circ \phi_m$ is harmonic for all $m \in \M \setminus D$, which means that
$\tl{u} \in \hop$.  Thus (c) implies~(d).

Now suppose (d) holds, so $\tl{u} \in \hop$.  Let $m \in \M \setminus D$.  Thus
$\tl{u} \circ \phi_m$ is a harmonic function and hence is equal to its Berezin
transform.  In other words,
\[
(\tl{u} \circ \phi_m)\til = \tl{u} \circ \phi_m = (u \circ \phi_m)\til,
\]
where the second equality comes from Lemma~\ref{32}(a).  Because the Berezin
transform is one-to-one (as we showed after the proof of Lemma~\ref{29}), the equation
above implies that
\[
\tl{u} \circ \phi_m = u \circ \phi_m.
\]
Evaluating both sides of this equation at $0$ shows that $\tl{u}(m) = u(m)$.  Thus
(d) implies~(b).  At this point in the proof we have shown that (a), (b), (c), and (d)
are equivalent.

A result of McDonald and Sundberg (\bibref{McS}, Proposition~5) states that for a
function $v \in \U$, the Toeplitz operator $T_v$ is compact if and only if
$\lim_{z \to \partial D} v(z) = 0$.  Applying this result with $v = \tl{u} - u$ shows
that (a) and (e) are equivalent.  Thus we now know that (a) through (e) are equivalent.

By Lemma~\ref{33}, eight times the integral in (f) equals
$\bigl(\Delta (u \circ \phi_z)\til\bigr)(0)$, which by Lemma~\ref{32}(a) equals
$\bigl(\Delta (\tl{u} \circ \phi_z)\bigr)(0)$, which by a standard calculation equals
$(1 - |z|^2)^2 (\Delta \tl{u})(z)$.  Thus (f) and (g) are equivalent.

Now suppose that (d) holds, so $\tl{u} \in \hop$.  Thus
\[
\bigl(\Delta (\tl{u} \circ \phi_m)\bigr)(0) = 0
\]
for every $m \in \M \setminus D$.  By Lemma~\ref{32}(c) (with $w=0$), this gives (g). 
Thus (d) implies~(g).

Now suppose (g) holds.  Fix $w \in D$ and $m \in \M \setminus D$.  Note that
$|\phi_z(w)| \to 1$ as $z \to m$ (recall that $w$ is fixed).  From
Lemma~\ref{32}(c) and our hypothesis (g) we now conclude that
$\bigl(\Delta (\tl{u} \circ \phi_m)\bigr)(w) = 0$.  Thus $\tl{u} \circ \phi_m$ is
harmonic on $D$.  Because $m$ was an arbitrary element of $\M \setminus D$, this
means that $\tl{u} \in \hop$.  Thus (g) implies (d), completing the proof that (a)
through (g) are equivalent.

To deal with (h), now suppose that $u$ is a finite sum of functions of the form
$u_1 \dots u_n$, where each $u_j$ is a bounded harmonic function on $D$.  For each
such $u_j$, the function $u_j \circ \phi_z$ is harmonic on $D$ for every $z \in D$. 
If $m \in \M$, then $u_j \circ \phi_z$ converges pointwise on $D$ to
$u_j \circ \phi_m$ as $z \to m$.  A pointwise convergent net of uniformly bounded
harmonic functions has the property that every partial derivative (of arbitrary order)
also converges pointwise to the appropriate partial derivative of the limit function. 
Applying this (and the appropriate product rule for partial derivatives) to $u$ gives
\begin{align*}
\bigl(\Delta (u \circ \phi_m)\bigr)(w) &=
\lim_{z \to m} \bigl(\Delta (u \circ \phi_z)\bigr)(w) \\
&= \lim_{z \to m} \frac{(1 - |\phi_z(w)|^2)^2 (\Delta u)(\phi_z(w))}{(1-|w|^2)^2},
\end{align*}
for every $w \in D, m \in \M$, where the second equality comes from a standard
calculation.  To prove that (c) is equivalent to~(h), now follow the pattern of the
proof showing that (d) is equivalent to~(g), using the last equality above in place of
Lemma~\ref{32}(c).
\qed

A continuous bounded function on $D$ can have Berezin transform in $C(\bar{D})$ without
itself being in $C(\bar{D})$ (of course, to say that a continuous function on $D$ is in
$C(\bar{D})$ means that it extends continuously to a function on $\bar{D}$).  To
construct an example, consider a continuous function $v$ on $[0, 1)$ that equals
$0$ most of the time (enough so that the average value of $v$ on the interval $[r, 1)$
tends to $0$ as $r$ increases to~$1$), but whose graph occasionally has a small bump
with height~$1$ (so that $v$ does not extend continuously to $[0,1]$).  Define a radial
function $u$ on $D$ by $u(z) = v(|z|)$. Then $\tl{u}$ extends continuously to $\bar{D}$
even though $u$ does not have this property.  The following corollary shows that
functions in $\U$ cannot behave in this fashion.

\begin{corollary} \label{80}
Suppose\/ $u \in \U$.  If\/ $\tl{u} \in C(\bar{D})$, then\/ $u \in C(\bar{D})$.
\end{corollary}

\proof
Suppose $\tl{u} \in C(\bar{D})$.  Then $\tl{u} \in \hop$ (because $u \circ \phi_m$
is constant on $D$ for every $m \in \M \setminus D$).  Because condition~(d) in
Theorem~\ref{36} holds, condition~(a) in the same theorem also
holds.  Condition~(a) and the continuity of $\tl{u}$ on $\bar{D}$ imply that
$u$ extends continuously to $\bar{D}$, completing the proof.
\qed

Suppose $u \in \U$.  In proving Theorem~\ref{36}, we used the McDonald-Sundberg
theorem that $T_u$ is compact if and only if
$u(z) \to 0$ as $z \to \partial D$ (\bibref{McS}, Proposition~5).  To provide an
easy proof of this theorem using our tools, note that Theorem~2.2 of
our paper \bibref{AxZ} asserts that $T_u$ is compact if and only if $\tl{u} \to 0$ as
$z \to \partial D$.  The equivalence of (a) and (d) in Theorem~\ref{36} shows that this
happens if and only if $u(z) \to 0$ as $z \to \partial D$, completing our proof of the
McDonald-Sundberg theorem.  (This is not a circular proof of the McDonald-Sundberg
theorem, as that result was used in the proof of Theorem~\ref{36} only in showing
that (e) is equivalent to (a); this equivalence is not used in the proof we have just
given).

The McDonald-Sundberg theorem proved in the paragraph above gives another
example of how $\U$ provides a more natural context than $L^\infty(D, dA)$ for many
Toeplitz operator questions.  Specifically, the McDonald-Sundberg theorem just proved
becomes false if the hypothesis that $u \in \U$ is weakened to the hypothesis
that $u \in L^\infty(D, dA)$---Sarason constructed an example, presented in Section~5
of~\bibref{Str}, of a function $u \in L^\infty(D, dA)$ such that $T_u$ is compact but
$|u(z)| = 1$ for all $z \in D$.

In Corollary~\ref{66}, we showed that $S - T_{\tl{S}}$ is in the commutator ideal
$\CT$ for every $S \in \T$.  This raises the question of when $S - T_{\tl{S}}$ is a
compact operator.  In the corollary below, we answer this question for $S$ lying in a
dense subset of~$\T$.  We do not know whether the hypothesis on $S$ in the corollary
below could be replaced by the weaker hypothesis that $S \in \T$.

\begin{corollary} \label{44}
Suppose\/ $S$ is a finite sum of operators of the form\/ $T_{u_1} \dots T_{u_n}$, where
each\/ $u_j \in \U$.  Then the following are equivalent:
\begin{subtheorem}
\item
$S - T_{\tl{S}}$ is a compact operator;

\item
$\tl{S} \in \hop$;

\item
${\displaystyle\lim_{z \to \partial D} (1 - |z|^2)^2 (\Delta \tl{S})(z) = 0}$.
\end{subtheorem}
\end{corollary}

\proof
In Theorem~2.2 of \bibref{AxZ}, we showed that a finite sum of finite products of
Toeplitz operators is compact if and only if its Berezin transform has limit $0$
on~$\partial D$.  Applying this result to the operator $S - T_{\tl{S}}$, whose
Berezin transform equals $\tl{S} - \Tilde{\Tilde{S}}$, we conclude that
$S - T_{\tl{S}}$ is compact if and only if
\begin{equation} \label{35}
\lim_{z \to \partial D} \tl{S}(z) - \Tilde{\Tilde{S}}(z) = 0.
\end{equation}%
The equivalence of conditions (a) and (c) in Theorem~\ref{36} (with $u = \tl{S}$) shows
that \eqref{35} holds if and only if $\tl{S} \in \hop$.  In other words, conditions
(a) and (b) above are equivalent.

Now suppose that (b) holds, so $\tl{S} \in \hop$.  Thus $\bigl(\Delta (\tl{S} \circ
\phi_m)\bigr)(0) = 0$ for every $m \in \M \setminus D$.  Proposition~\ref{31}(c) (with
$w=0$) now tells us that
$\lim_{z \to \partial D} (1 - |z|^2)^2 (\Delta \tl{S})(z) = 0$.  In other
words, (b) implies (c).

Now suppose that (c) holds.  Fix $w \in D$ and $m \in \M \setminus D$.  Note that
$|\phi_z(w)| \to 1$ as $z \to m$ (recall that $w$ is fixed).  From
Proposition~\ref{31}(c) and our hypothesis (c) we now conclude that
$\bigl(\Delta (\tl{S} \circ \phi_m)\bigr)(w) = 0$.  Thus $\tl{S} \circ \phi_m$ is
harmonic on $D$.  Because $m$ was an arbitrary element of $\M \setminus D$, this
means that $\tl{S} \in \hop$.  Thus (c) implies (b), completing the proof.
\qed

\section{Asymptotic Multiplicativity}

The Berezin transform is not multiplicative even over the space of harmonic functions.
However, $\wt{uv}(z)-\tl{u}(z)\tl{v}(z)\to 0$ as $z \to \partial D$ for some pairs of
functions $u, v$.  In this section we describe when this happens for bounded
harmonic functions.  Note that if $u$ and $v$ are harmonic, then $\tl{u} = u$ and
$\tl{v} = v$, so we want to know when $\wt{uv}$ is approximately equal to $uv$
near~$\partial D$.

Our key tool in proving Theorem~\ref{45} will be Theorem~\ref{36}.  However we will
also need the following two lemmas.

\begin{lemma} \label{40}
If\/ $u, v$ are bounded and harmonic on $D$, then
\[
\wt{uv}(z)-u(z)v(z)=
({H_{\bar{u}}}^*H_{v}+{H_{\bar{v}}}^*H_{u})\til(z)
\]
for every\/ $z \in D$.
\end{lemma}

\proof
Suppose $u$ and $v$ are bounded and harmonic on $D$.  There there are four functions
$f_{1}, f_{2}, g_{1}, g_{2} \in L^2_a(D)$ such that
$u = f_{1}+\bar{f}_{2}$ and $v = g_{1}+\bar{g}_{2}$. Let $z \in D$.  Then
\begin{align*}
({H_{\bar{u}}}^* H_v)\til(z) &= \ip{{H_{\bar{u}}}^* H_v k_z, k_z} \\
&= \ip{H_v k_z, H_{\bar{u}} k_z } \\
&= \ip{ (1-P)\bigl((g_1 + \overline{g_2})k_z\bigr),
(1-P)\bigl((\overline{f_1} + f_2)k_z\bigr) } \\
&= \ip{ (1-P)(\overline{g_2}k_z), (1-P)(\overline{f_1}k_z) } \\
&= \ip{ \overline{g_2}k_z - \overline{g_2(z)}k_z,
\overline{f_1}k_z - \overline{f_1(z)}k_z } \\
&= \ip{ f_1 \overline{g_2} k_z, k_z } - f_1(z) \overline{g_2(z)}.
\end{align*}
{\parfillskip=0pt plus1fil
A similar formula holds for $({H_{\bar{v}}}^* H_u)\til(z)$.  Adding these two formulas
gives
\begin{align*}
({H_{\bar{u}}}^*H_{v} + {H_{\bar{v}}}^* H_{u})\til(z)
&= \ip{ ( f_1 \overline{g_2} + \overline{f_2} g_1 ) k_z, k_z}
- f_1(z) \overline{g_2(z)} - \overline{f_2(z)} g_1(z) \\
&= \ip{ (f_1 + \overline{f_2}) (g_1 + \overline{g_2}) k_z, k_z } \\
& \quad\quad {}- (f_1(z) + \overline{f_2(z)})(g_1(z) + \overline{g_2(z)}) \\
&= \ip{ uv k_z, k_z} - u(z) v(z) \\
&= \wt{uv}(z)-u(z)v(z),
\end{align*}
as desired.
\qed

}Although the following lemma is probably well known, we were unable to locate a proof
in the literature.  Thus we have included a proof.

\begin{lemma} \label{41}
Suppose\/ $u, v$ are harmonic on\/ $D$.  Then\/ $uv$ is harmonic on\/ $D$ if and only
if at least one of the following conditions holds:

\begin{subtheorem}
\item
$u$ and\/ $v$ are both analytic on\/ $D$;  
  
\item
$\bar{u}$ and\/ $\bar{v}$ are both analytic on\/ $D$;

\item
there exist complex numbers\/ $\alpha, \beta$, not both\/ $0$, such that\/
$\alpha u +\beta v$ and\/ $\bar{\alpha} \bar{u} - \bar{\beta} \bar{v}$ are both
analytic on\/ $D$.
\end{subtheorem}
\end{lemma}

\proof
Because $u$ and $v$ are harmonic on $D$, an elementary computation shows that 
\[
\Delta(uv) = 4\left(
\frac{\partial u}{\partial \bar{z}} \frac{\partial v}{\partial z} +
\frac{\partial u}{\partial z} \frac{\partial v}{\partial \bar{z}} 
\right)
\]
on $D$.  Thus $uv$ is harmonic if and only if
\begin{equation} \label{38}
\frac{\partial u}{\partial \bar{z}}
\frac{\partial v}{\partial z} = -
\frac{\partial u}{\partial z} \frac{\partial v}{\partial \bar{z}}.
\end{equation}%

Clearly (a) implies that $uv$ is harmonic, as does
(b).  Condition (c) can be restated to say that there exist complex numbers $\alpha,
\beta$, not both $0$, such that
\begin{equation} \label{39}
\alpha \frac{\partial u}{\partial \bar{z}} =
- \beta \frac{\partial v}{\partial \bar{z}} \quad\text{and} \quad
\alpha \frac{\partial u}{\partial z} = \beta \frac{\partial v}{\partial z}.
\end{equation}%
Thus (c) implies \eqref{38}, proving one direction of the lemma.

To prove the other direction, we use an argument from the proof of Theorem~1
of~\bibref{AxC}.  Suppose that $uv$ is harmonic, so
\eqref{38} holds.  Because
$u$ and $v$ are harmonic, $\frac{\partial u}{\partial z}$ and 
$\frac{\partial v}{\partial z}$ are analytic; furthermore,
$\frac{\partial u}{\partial\bar{z}}$ and $\frac{\partial v}{\partial\bar{z}}$ are
conjugate analytic.
Let
\[
\Omega =\{w\in D: \frac{\partial v}{\partial z}(w) \neq 0 \text{\ and\ }
\frac{\partial v}{\partial \bar{z}}(w) \neq 0 \}.
\]

First consider the case where $\Omega$ is the empty set.  Then either $v$ is analytic
or $v$ is conjugate analytic.  If $v$ is analytic, then \ref{38} implies
$\frac{\partial u}{\partial \bar{z}} \frac{\partial v}{\partial z} = 0$, which
implies that either $u$ is analytic (so (a) holds) or $v$ is conjugate
analytic (so (c) holds with $\alpha = 0$ and $\beta = 1$).  Similarly, if $v$ is
conjugate analytic, then either (b) or (c) holds, completing the proof when $\Omega$
is the empty set.

Now suppose $\Omega$ is not the empty set.  Then $\Omega$ is a dense open subset
of~$D$.  On $\Omega$, we can rewrite \ref{38} as
\[
- \frac{\frac{\partial u}{\partial \bar{z}} }{\frac{\partial v}{\partial \bar{z}}}
= \frac{\frac{\partial u}{\partial z}}{\frac{\partial v}{\partial z}}.
\]
The left side of this  equation is a conjugate analytic function on~$\Omega$. The right
side is an analytic function on~$\Omega$. Thus both sides equal the same constant
function on~$\Omega$.  We conclude that for some constant $\beta$, we must have 
$\frac{\partial  u}{\partial \bar{z}}
= -\beta \frac{\partial v}{\partial \bar{z}}$
and $\frac{\partial u}{\partial z} = \beta         
\frac{\partial v}{\partial z}$.  Thus \eqref{39} holds (with $\alpha = 1$) and hence
(c) holds, completing the proof.
\qed

Now we can describe when the Berezin transform is asymptotically multiplicative on
harmonic functions.  For partial results on when the Berezin transform is
multiplicative on $\B(L^2_a)$, see Kili\c{c}'s paper~\bibref{Kil}.

\begin{theorem} \label{45}
Suppose\/ $u$ and\/ $v$ are bounded and harmonic on\/ $D$.
Then the following conditions are equivalent:

\begin{subtheorem}

\item
${\displaystyle \lim_{z \to \partial D} \wt{uv}(z) - u(z) v(z)= 0}$;

\item
${\displaystyle \lim_{z \to \partial D} (1-|z|^2)^{2} \Delta(uv) (z)=0}$;

\item
$uv \in \hop$;

\item
$2T_{uv} - T_u T_v - T_v T_u$ is compact;

\item
for each\/ $m \in \M \setminus D$, at least one of the following conditions holds:
\begin{subsubtheorem}
\item
$u \circ \phi_m$ and\/ $v \circ \phi_m$ are both in\/ $H^{\infty}$;
  
\item
$\bar{u} \circ \phi_m$ and\/ $\bar{v} \circ \phi_m$ are both in\/ $H^{\infty}$;

\item
there exist complex numbers\/ $\alpha, \beta$, not both\/ $0$, such that\\
$\alpha u \circ \phi_m +\beta v \circ \phi_m$ and\/
$\bar{\alpha} \bar{u}\circ \phi_m - \bar{\beta} \bar{v} \circ \phi_m$ are both in\/
$H^{\infty}$.
\end{subsubtheorem}
\end{subtheorem}
\end{theorem}

\proof
The equivalence of (a), (b), and (c) follows from the equivalence of conditions (a),
(h), and (c) in Theorem~\ref{36}.

The equivalence of (a) and (d) follows from Lemma~\ref{40} and Theorem~2.2 of
\bibref{AxZ} along with the identity $2T_{uv} - T_u T_v - T_v T_u =
{H_{\bar{u}}}^*H_{v} + {H_{\bar{v}}}^* H_{u}$.

Finally, the equivalence of (c) and (e) follows from Lemma~\ref{41}.
\qed

As an application of the theorem above, we now show how it can be used to
give an easy proof of the characterization of the functions $f, g \in H^\infty$ such
that $T_{\bar{f}} T_g - T_g T_{\bar{f}}$ is compact.  Suppose $f, g \in H^\infty$. 
Then
\[
2 T_{\bar{f}g} - T_{\bar{f}} T_g - T_g T_{\bar{f}} =
T_{\bar{f}} T_g - T_g T_{\bar{f}}.
\]
Thus by Theorem~\ref{45} (with $u = \bar{f}$ and $v = g$),
$T_{\bar{f}} T_g - T_g T_{\bar{f}}$ is compact if and only if
\[
\lim_{z \to \partial D} (1-|z|^2)^{2} \Delta(\bar{f}g)(z) = 0.
\]
Because $\Delta = 4(\partial/\partial z)(\partial/\partial \bar{z})$, we see
that $T_{\bar{f}} T_g - T_g T_{\bar{f}}$ is compact if and only if
\begin{equation} \label{81}
\lim_{z \to \partial D} (1-|z|^2)^{2} f'(z) g'(z) = 0.
\end{equation}%
This result was originally proved by Zheng \bibref{Zh2} using other methods; also see
\bibref{AxG} for additional conditions that are equivalent to~\eqref{81}.

\newpage

\section*{References}

\begin{list}{\arabic{referencec}.\hfill}{\usecounter{referencec} 
\setlength{\topsep}{12pt plus 3pt minus 3pt}
\setlength{\partopsep}{0pt}
\setlength{\labelwidth}{1.5\parindent}
\setlength{\labelsep}{0pt}
\setlength{\leftmargin}{1.5\parindent}
\setlength{\parsep}{1ex} }

\raggedright

\item \label{AFR}
Patrick Ahern, Manuel Flores, and Walter Rudin,
An invariant volume-mean-value property,
\textsl{J. Funct. Anal.} 111 (1993), 380--397.

\item \label{Ax2}
Sheldon Axler,
Bergman spaces and their operators,
\textsl{Surveys of Some Recent Results in Operator Theory}, vol.~1, edited by John
B.~Conway and Bernard B.~Morrel, Pitman Research Notes in Mathematics, 1988, 1--50.

\item \label{AxC}
Sheldon Axler and \u{Z}eljko \u{C}u\u{c}kovi\'c,
Commuting Toeplitz operators with harmonic symbols,
\textsl{Integral Equations Operator Theory} 14 (1991), 1--12.

\item \label{AxG}
Sheldon Axler and Pamela Gorkin,
Algebras on the disk and doubly commuting multiplication operators,
\textsl{Trans. Amer. Math. Soc.} 309 (1988), 711--723.

\item \label{ASh}
Sheldon Axler and Allen Shields,
Extensions of analytic and harmonic functions,
\textsl{Pacific J. Math.} 145 (1990), 1--15.

\item \label{AxZ}
Sheldon Axler and Dechao Zheng,
Compact operators via the Berezin transform,
preprint.

\item \label{Eng}
Miroslav Engli\v s,
Functions invariant under the Berezin transform,
\textsl{J. Funct. Anal.} 121 (1994), 233--254.

\item \label{Gar}
John B. Garnett,
\textsl{Bounded Analytic Functions},
Academic Press, 1981.

\item \label{Hof}
Kenneth Hoffman,
Bounded analytic functions and Gleason parts,
\textsl{Ann. of Math.} 86 (1967), 74--111.

\item \label{Kil}
Semra Kili\c{c},
The Berezin symbol and multipliers of functional Hilbert spaces,
\textsl{Proc. Amer. Math. Soc.} 123 (1995), 3687--3691.

\item \label{McS}
G. McDonald and C. Sundberg,
Toeplitz operators on the disc,
\textsl{Indiana Univ. Math. J.}  28 (1979), 595--611.

\item \label{Str}
Karel Stroethoff,
Compact Toeplitz operators on Bergman spaces,
preprint.

\item \label{Sun}
Carl Sundberg,
Exact sequences for generalized Toeplitz operators,
\textsl{Proc. Amer. Math. Soc.} 101 (1987), 634--636.

\item \label{Zh2}
Dechao Zheng,
Hankel operators and Toeplitz operators on the Bergman space,
\textsl{J. Funct. Anal.} 83 (1989), 98--120.

\end{list}

\bigskip

\bigskip

\pagebreak[3]
\noindent
{\scshape
Sheldon Axler \\
Department of Mathematics \\
San Francisco State University \\
San Francisco, CA 94132 USA

\bigskip

\noindent
Dechao Zheng \\
Department of Mathematics \\
Vanderbilt University \\
Nashville, TN 37240 USA
}

\bigskip

\noindent
\textsl{e-mail}: \texttt{axler@math.sfsu.edu} and \texttt{zheng@math.vanderbilt.edu}

\bigskip

\noindent
\textsl{Axler www home page}: \texttt{http://math.sfsu.edu/axler}

\end{document}